\documentclass{elsarticle-modified}
\journal{Appl. Math. Comput.}

\usepackage[english]{babel}
\usepackage{amsmath,amssymb,amsthm}
\usepackage{xcolor}

\def\nE         {\ensuremath{{\mathcal E}}}
\def\nS         {\ensuremath{{\mathcal S}}}
\def\nL         {\ensuremath{{\mathcal L}}}

\newcommand{\nP}{\mathcal{P}}

\newcommand\tnS{\tilde{\mathcal{S}}}
\newcommand\tnL{\tilde{\mathcal{L}}}
\newcommand\Order{{\mathcal{O}}}

\newcommand{\RR}{{{\mathbb{R}}\vphantom{|}}}

\newcommand{\ee}{{\mathrm e}}

\begin{document}

\begin{frontmatter}

\title{Construction of adaptive exponential multi-operator splitting methods}

\author[wpi,mmm]{O. Koch}\corref{mycorrespondingauthor}
\ead{othmar@othmar-koch.org}

\author[tu]{K. Acar}
\ead{e11813035@student.tuwien.ac.at}

\author[tu]{W. Auzinger}
\ead{winfried.auzinger@tuwien.ac.at}

\author[wpi,fhtw,mmm]{F. Kupka}
\ead{kupka@technikum-wien.at}

\author[tu]{D. Hoffmann}
\ead{daniel.hoffmann@tuwien.ac.at}

\author[tu]{B. Moser}
\ead{e1527052@student.tuwien.ac.at}

\address[wpi]{Wolfgang Pauli Institute, Oskar-Morgenstern-Platz 1, A--1090 Wien, Austria}

\address[tu]{Institut f{\"u}r Analysis und Scientific Computing (E101), TU Wien, Wiedner Hauptstra{\ss}e 8--10, A--1040 Wien, Austria}

\address[fhtw]{Dept. Applied Mathematics and Physics, Univ. of Applied Sciences Technikum Wien, H{\"o}chst{\"a}dtplatz 6, A--1200 Wien, Austria}

\address[mmm]{Research platform ``MMM -- Mathematics, Magnetism, Materials'', Faculty of Mathematics, Univ. Vienna, Oskar-Morgenstern-Platz 1, A--1090 Wien, Austria}

\cortext[mycorrespondingauthor]{Corresponding author}

\begin{abstract}
We construct splitting methods suitable for the solution of the equations of magnetohydrodynamics (MHD).
Due to the physical significance of the involved operators, splittings into three or even four operators
with positive coefficients are appropriate for a physically correct and efficient
solution of the equations. To efficiently obtain an accurate solution approximation, adaptive
choice of the time-steps is important particularly in the light of the unsmooth
dynamics of the system. Thus, we construct new method coefficients in conjunction
with associated error estimators by optimizing the leading local error term.
As a proof of concept, we demonstrate that adaptive splitting faithfully reflects
the solution behavior also in the presence of a shock-like behavior for the viscous Burgers
equation, which serves as a simplified model problem displaying several features of
the Navier--Stokes equation for incompressible flow.
\end{abstract}


\end{frontmatter}

\section{Introduction}\label{sec:intro}

We construct and test operator splitting methods suitable for the equations of magnetohydrodynamics (MHD)
(\ref{mhdeqns1}) in astrophysics, see \cite{galtier16b}.
For the last few decades astrophysical research has experienced a unique phase of growth.
This has been driven by the collection of ever more accurate and more complete experimental data on individual
and also large classes of astrophysical objects (e.g., by the Hubble Space Telescope, ESO's VLT, etc.)
as well as on the development of ever more sophisticated computational techniques to solve dynamical
equations describing their physical properties. Many astrophysical objects can be described as a single or
multicomponent fluid, typically in the state of a plasma. This is also the case for stars such as our Sun,
which are in the focus of our interest. In most astro- and geophysical
problems the dynamics of a fluid is described by the conservation laws of hydrodynamics.
A complete set of equations \cite{hujran01,Muthsam2010} includes gas density, concentration of the second
species, flow velocity, and total energy density. Where forces induced by an external magnetic field or one
generated by the flow of the fluid cannot be neglected, a dynamical equation is added to model the interaction
between the magnetic field and the flow. The associated equations of MHD (\ref{mhdeqns1})
introduce the constraint of divergence-freeness of the magnetic field (see also \cite{brun17,stein12} which
provide an overview on research results from MHD simulations and experimental studies of our Sun and other stars).

The applications we are interested in are simulations of the solar surface. This field of research is relevant to the
understanding of astrophysical phenomena such as solar activity and thus further improves our knowledge on
solar forcing which is also of interest to climate modeling.

The implementation our developments are going to be integrated into in the future is the ANTARes code
\cite{Muthsam2010}. This is a software tool to solve, among others,
the Euler or alternatively the Navier-Stokes equations of a single or two mixing species,
augmented by equations for radiative transport and an equation of state.
This is currently being extended to the equations of MHD and will
be equipped with adaptive splitting methods to increase the performance.

\section{Motivation --- Magnetohydrodynamics}\label{sec:mhd}

A full set of the equations of MHD reads
\begin{subequations}\label{mhdeqns1}
{\footnotesize
\begin{eqnarray}\label{hydro1}
\underbrace{\partial_t \left( \begin{array}{c} \rho \\
                                                              \rho c \\
                                                              \rho \boldsymbol{u} \\
                                                              e \\
                                        \end{array}\right)}_{\partial_t y(t)} &=&
   \underbrace{ - \nabla \cdot \left( \left( \begin{array}{c} \rho \boldsymbol{u}  \\
                                                                    \rho c \boldsymbol{u}  \\
                                                                    \rho \boldsymbol{u}  \otimes \boldsymbol{u}  \\
                                                                    e \boldsymbol{u}
                                                        \end{array}\right)
                                              +\left( \begin{array}{c} 0 \\
                                                                           0 \\
                                                                           P \mathbf{Id}\\
                                                                           P \boldsymbol{u}
                                                        \end{array}\right)
                                              -\left( \begin{array}{c} 0 \\
                                                                             0 \\
                                                                             \sigma \\
                                                                             \boldsymbol{u}  \cdot \sigma
                                                        \end{array}\right)\right)
                                              +\left( \begin{array}{c} 0 \\
                                                                   0 \\
                                                                   \rho \boldsymbol{g} \\
                                                                   \rho \boldsymbol{g}\cdot \boldsymbol{u}
                                                       \end{array}\right) }_{A(y(t))=(A_1+A_2+A_3+A_4)(y(t))} \nonumber\\
&& +  \underbrace{\nabla \cdot \left( \begin{array}{c} 0 \\
                                                                       \rho \kappa_c \nabla c \\
                                                                       0 \\
                                                                       K \nabla T
                                               \end{array} \right)}_{B(y(t))}\,
+  \underbrace{\left(\begin{array}{c} 0 \\ 0 \\ \boldsymbol{j}\times\boldsymbol{B} \\
\boldsymbol{u}\cdot(\boldsymbol{j}\times\boldsymbol{B})+\eta \boldsymbol{j}^2\end{array}\right)}_{C(y(t))},
\end{eqnarray}
}
coupled to the magnetic induction equation
\begin{eqnarray}\label{magnetic1}
\partial_t \boldsymbol{B} = \mathrm{curl}(\boldsymbol{u}\times\boldsymbol{B}) + \eta \Delta \boldsymbol{B}.
\end{eqnarray}
\end{subequations}
Here, the symbols for variables and parameters have their standard meaning in hydrodynamics.
The dynamical variables represent the following physical quantities: $\mathbf{u}$ is the flow velocity,
$\rho$ is the gas density, $c$ describes the concentration of the second species, and $e$ represents
the total energy density (sum of kinetic and internal energy to which the magnetic energy could be added separately, see
also \cite{stein12}).

$\boldsymbol{B}\in\RR^3$ is the magnetic field and in non-ideal MHD the magnetic
diffusivity $\eta \neq 0$. Note that $\boldsymbol{B}$ has to satisfy the constraint
\begin{equation}\label{div0}
      \mathrm{div}\boldsymbol{B}=0 \qquad \forall\, t.
\end{equation}
In numerical computations, this \emph{divergence-freeness} of the magnetic field should
be exactly or at least very accurately preserved, see \cite{fuchsetal09b,toth00}.

The structure of the model equations naturally suggests an additive splitting of (\ref{mhdeqns1}) consisting of terms each
characterized by different mathematical properties: the nonlinear, hyperbolic advection
terms ($\rho \boldsymbol{u}, \rho c \boldsymbol{u}, \rho \boldsymbol{u}  \otimes \boldsymbol{u},
e \boldsymbol{u}$); the pressure-containing terms which often induce time-step
limits for explicit time integration methods due to the short time-scale associated with
sound waves that are introduced through these terms; the source terms containing
the gravitational acceleration $\boldsymbol{g}$; the parabolic diffusion terms;
and terms related to the magnetic induction. One natural way to split can thus be defined by the suboperators
$\mathcal{A}=(A_1,\mathrm{curl}(\boldsymbol{u}\times \boldsymbol{B}))^T$, $\mathcal{B}=(A_2+A_3+A_4,0,0,0)^T$,
$\mathcal{C}=(B(y(t)), \eta \Delta \boldsymbol{B})^T$, $\mathcal{D}=(C(y(t),0,0,0)^T$.
Furthermore, it could be considered to additionally split off the terms
containing the magnetic induction in $\mathcal{C}$ and propagate these simultaneously with
$\mathrm{curl}(\boldsymbol{u}\times \boldsymbol{B})$,
or to treat the operator $C$ together with $\mathcal{B}$. In this way, splitting into
3--7 operators can be considered.
Also, terms could be grouped differently based on physical considerations
and employed discretization.
The decision which splitting is most favorable depends also on details of the implementation
in the ANTARes code and cannot easily be made a priori. Thus, the present paper presents
splittings into three and four operators, respectively, which can be employed in a future
implementation for real world problems.

\section{Splitting methods}\label{sec:splitting}

We will apply additive splitting methods for the time integration.
Consider a nonlinear partitioned evolution equation
\begin{equation}\label{partitioned}
\partial_t u(t) = A_1(u(t)) + \cdots + A_n(u(t)),\qquad u(0)=u_0 \quad \mbox{given}.
\end{equation}
A general splitting into the $n$ operators is defined as follows (see also \cite{haireretal02b}): Denote the
flow generated by the operator $A_j$ by $\nE_j$, i.e.
$u_j(t)=\nE_j(t,u_j^0)$ solves the abstract initial value problem
%
$$ \partial_t u_j(t) = A_j(u_j(t)),\quad u_j(0)=u_j^0,\quad j=1,\dots,n.$$
Then, a splitting with $k$ compositions, $\nS(h,u_0): u_0 \mapsto u_1 \approx u(h)$ is defined by
\begin{eqnarray}
u_1 &=& \prod_{\nu=1}^{k}\prod_{\ell=1}^{n} \nE_{\ell}(a_{\ell,\nu}h,\cdot) u_0 \label{eq:splitting}\\
    &=& \nE_{n}(a_{n,k}h,\cdot)\circ \cdots \circ \nE_1(a_{1,k}h,\cdot)\circ \cdots
        \circ \nE_{n}(a_{n,1}h,\cdot) \circ \cdots \circ \nE_1(a_{1,1}h,\cdot) u_0.\nonumber
\end{eqnarray}
Here, $\circ$ stands for the composition of mappings.
The coefficients $a_{\ell,\nu}$ are suitably determined from the \emph{order conditions}
(a system of polynomial equations \cite{haireretal02b}) such that the approximation
has a desired order. To set up the order conditions is an issue of its own,
based on theoretical results about Lie algebrae \cite{auzingeretal16b}, automatic
tools have been implemented in Julia and in Maple to realize this
efficiently (without missing any relationships and without redundancy),
see~\cite{maple16,moskau19}.

An important aspect of determining the coefficients is also to optimize
the relationship between computational effort and numerical error
(in the \emph{leading local error term} $C\tau^{p+1}$, $C$ should be
as small as possible)\footnote{In fact, $C$ strongly depends on the
considered problem, it is a linear combination of commutators of the involved
operators. Thus, a globally optimal method cannot be determined, but
a generic quantity (Euclidean norm of the coefficients in the leading local error term) has been defined in
\cite{part1}, see also \cite{splithp}.}. Extensive optimizations have been conducted earlier,
yielding optimized splittings into two ($n=2$) and three ($n=3$) operators,
see \cite{part1} for explanations and \cite{splithp} for a collection
of coefficients which is continually extended. Splittings into a higher number of operators
has also been discussed recently in \cite{ostmeyer23}, where it is shown how two-operator
splittings can be extended to an arbitrary even number of operators and the performance
of the resulting schemes is assessed.

Another important feature in the considered model is the fact that operators
which are either parabolic or hyperbolic, must not be propagated backward in time,
implying the need for positive real coefficients\footnote{For parabolic operators,
an alternative is to choose complex coefficients with positive real parts \cite{Hansen09}.
This cannot be used, however, for an implementation in ANTARes.}.

In the construction and application of splitting methods for MHD we will focus on
methods of order two, since a consistent discretization of the boundary
conditions is difficult for higher-order schemes and the solution's lack of
smoothness also suggests to limit the order.

\section{Adaptive time-stepping}\label{sec:adapt}

\subsection{Local error estimation}\label{errest}

In this section, we briefly describe three classes of local error estimators
which may serve as our basis for adaptive time-stepping and which have different
advantages depending on the context in which they are applied.
The error estimates described below
are \emph{asymptotically correct}, i.e., the deviation of the error estimator
from the true error tends to zero faster than does the error.

\emph{Embedded pairs of splitting formulae} of different orders have first been considered in
\cite{knth10b} and are based on reusing a number of evaluations from
the basic integrator.
Since this approach involves a second and a third order method, consistent
discretization at the boundary has to be ensured. This raises a number of issues
which need to be resolved and thus, this approach is not pursued at present.

For \emph{adjoint pairs of formulae} of odd order $p$ or Milne pairs, an asymptotically correct error
estimator can be computed at the same cost as for the basic method, see \cite{part1}.
Since a basic method of order 1 is required to obtain an asymptotically correct
error estimator, we do not pursue this approach at present because of
the expected loss in accuracy and efficiency.

\paragraph{Milne pairs} \label{milne} In the context of multistep methods for ODEs,
the so-called Milne device is a well-established technique for constructing
pairs of schemes. In our context, one may aim for finding a pair
$ (\nS,\tnS) $ of schemes of the same type, with equal $ s $ and $ p $, such that their local errors
$ \nL,\tnL $ are related according to
\begin{subequations} \label{eq:le-milne}
\begin{align}
\nL(h,u)  &= ~~C(u)\,h^{p+1} + \Order(h^{p+2}), \label{eq:le-milne1} \\
\tnL(h,u) &= \gamma\,C(u)\,h^{p+1} + \Order(h^{p+2}), \label{eq:le-milne2}
\end{align}
\end{subequations}
with $ \gamma \not= 1 $. Then, the additive scheme
\begin{equation*}
{\bar\nS}(h,u) = -\tfrac{\gamma}{1-\gamma}\,\nS(h,u) + \tfrac{1}{1-\gamma}\,\tnS(h,u)
\end{equation*}
is a method of order $ p+1 $, and
\begin{equation*}
\nS(h,u) - {\bar\nS}(h,u) = \tfrac{1}{1-\gamma}\,\big( \nS(h,u) - \tnS(h,u) \big)
\end{equation*}
provides an asymptotically correct local error estimate for $ \nS(h,u) $.
This construction principle requires detailed knowledge on the error structure
and is thus limited to low order methods.

In the construction, we additionally desire some level of embedding of coefficients
(see \cite{knth10b}) such that some of the subflows need only be computed
once and can be used in both methods.

\subsection{Step-size selection}\label{adapt}

Based on a local error estimator, the step-size is adapted such that the tolerance
is expected to be satisfied in the following step. If $h_{\text{old}}$ denotes the
current step-size, the next step-size $h_{\text{new}}$ in an order $p$ method is predicted as (see \cite{haireretal87,recipes88})
\begin{equation}
\label{step-selct}
h_{\text{new}} = h_{\text{old}} \cdot \min\Big\{\alpha_{\text{max}},\max\Big\{\alpha_{\text{min}},
\alpha\,\Big(\dfrac{\text{tol}}{\nP(h_{\text{old}})}\Big)^{\frac{1}{p+1}}\,\Big\}\Big\},
\end{equation}
where, similarly as in ODE solvers, we may choose for example
$ \alpha = 0.9 $, $ \alpha_{\text{min}} = 0.25 $, $ \alpha_{\text{max}} = 4.0 $,
and $\nP(h_{\text{old}})$ is an estimator for the local
error arising in the previous time-step.
This simple strategy incorporates safety factors to avoid an oscillating and unstable behavior.

\section{New methods and error estimators}\label{sec:construction}

In the following, we give some sets of coefficients which
have been determined so far. For splitting into four operators, optimized
methods have been determined under the nonnegativity condition.
For three operators, additionally a Milne pair for local error
estimation was found.

For readability, the coefficients are given with a truncated
mantissa, versions with double precision accuracy to be used for implementation are provided at \cite{splithp}.
Coefficients with a shorter mantissa are exact.
The coefficients have been obtained by Maple. Where no unique solution was defined,
minimization of a local error measure (LEM, defined as the Euclidean norm of
the coefficients in the residual with respect to the third-order condition,
see \cite{part1}) was conducted. In some cases, where the Maple routines could
not determine solutions of the nonlinear equations, the external package
\texttt{DirectSearch} provided at\\
\texttt{https://de.maplesoft.com/Applications/Detail.aspx?id=87637}\\
was employed.

\subsection{Splitting into four operators}\label{subsec:4ops}

\begin{itemize}
\item The coefficients of Strang splitting for four operators are given by
\begin{align}
  &a_{1,1} = 0, &&a_{1,2} = 0, &&a_{1,3} = 0, &&a_{1,4} = 1,\nonumber\\
  &a_{2,1} = 0, &&a_{2,2} = 0, &&a_{2,3} = 0.5, &&a_{2,4} = 0.5, \nonumber\\
  &a_{3,1} = 0, &&a_{3,2} = 0.5, &&a_{3,3} = 0, &&a_{3,4} = 0.5, \nonumber\\
  &a_{4,1} = 0.5, &&a_{4,2} = 0, &&a_{4,3} = 0, &&a_{4,4} = 0.5.\label{coeffs:strang}
\end{align}
To get an insight into the expected usefulness of the methods,
the \emph{local error measure (LEM)} introduced in \cite{part1} has been determined:
\begin{equation*}
  \text{LEM} = 2.6
\end{equation*}
\item Under the physically mandatory non-negativity constraint, the coefficients in (\ref{eq:splitting})
were determined by optimization and are given in Table~\ref{coeffs:optim4}. For any value
of $h$, this yields a second-order method, the local
error measure is optimized analytically for the value $h= \tfrac{169+3g-g^2}{12g}\approx 0.22633512,\ g=\sqrt[3]{72+\sqrt{4831993}}$.
\begin{table}
\begin{center}
\[\begin{tabular}{|c||c|c|c|c|}
\hline
 & $a_{1,*}$ & $a_{2,*}$ & $a_{3,*}$ & $a_{4,*}$\\
\hline
1 & 0 & $h$ & 0 & 0.5\\
\hline
2 & 0 & $0.5-h$ & 0.5 & 0\\
\hline
3 & 1 & 0 & 0.5 & $h$\\
\hline
4 & 0 & 0.5 & 0 & $0.5-h$\\
\hline
\end{tabular}\]
\caption{A four-operator splitting with positive coefficients (Method (I)).\label{coeffs:optim4}}
\end{center}
\end{table}

The analysis shows
\begin{equation*}
  \text{LEM} = 2.1,
\end{equation*}
The LEM is slightly smaller than for the Strang coefficients.

When we allow some additional freedom in the choice of $h$ to allow for the construction of
a Milne pair with Method (I) from Table~\ref{coeffs:optim4}, an associated Method (II) can be found
for the value $h\approx 0.3790984677886843$, where $\kappa=2.176315684585609$ such that $\kappa$((I)-(II)) is an asymptotically correct error estimator
for the basic method (I). The coefficients are given in Table~\ref{coeffs:A4332Milne}.
\begin{table}
\begin{center}
\begin{tabular}{|c||c|c|}
\hline
 & $a_{1,*}$ & $a_{2,*}$\\
\hline
1 & 0.13044731 & 0 \\
\hline
2 & 0 & 0.014157681\\
\hline
3 & 0 & 0.21691004\\
\hline
4 & 0 & 0.31230714\\
\hline
5 & 0.35103245 & 0\\
\hline
6 & 0.061895092 &  0.042865729\\
\hline
7 & 0 &  0.017373894\\
\hline
8 & 0.20166638 & 0\\
\hline
9 & 0 &  0.38686107\\
\hline
10 & 0 & 0.0095244307\\
\hline
\hline
 & $a_{3,*}$ & $a_{4,*}$\\
\hline
1 & 0.13044731 & 0.026400543\\
\hline
2 & 0 & 0.056956817\\
\hline
3 & 0 & 0.46001750\\
\hline
4 & 0.41292754 & 0\\
\hline
5 & 0 & 0\\
\hline
6 & 0 & 0\\
\hline
7 & 0.060239624 & 0\\
\hline
8 & 0.39638550 & 0.090841757\\
\hline
9 & 0 & 0.29331319\\
\hline
10 & 0 & 0.012230558\\
\hline
\end{tabular}
\caption{A four-operator splitting with positive coefficients (Method (II)), forming a Milne pair with the method from Table~\ref{coeffs:optim4}
(Method (I)).\label{coeffs:A4332Milne}}
\end{center}
\end{table}
%
\item A set of coefficients with five stages for splitting into four operators optimized under the
non-negativity constraint are given in Table~\ref{coeffs:optim5}.

\begin{table}
\begin{center}
\begin{tabular}{|c||c|c|}
\hline
 & $a_{1,*}$ & $a_{2,*}$\\
\hline
1 & 0.19859897 & 0.20567399\\
\hline
2 & 0.16188373 & 0.053687812\\
\hline
3 & 0.00000254 & 0.44666619\\
\hline
4 & 0.47832 & 0.094242\\
\hline
5 & 0.16119 & 0.19973\\
\hline
\hline
 & $a_{3,*}$ & $a_{4,*}$\\
\hline
1 & 0.15538119 & 0.43051849\\
\hline
2 & 0.43781080 & 0.071274504\\
\hline
3 & 0.13242 & 0.060827\\
\hline
4 & 0.067038 & 0.43738\\
\hline
5 & 0.20735 & 0\\
\hline
\end{tabular}
\caption{A four-operator splitting with positive coefficients.\label{coeffs:optim5}}
\end{center}
\end{table}

The resulting LEM is
\begin{equation*}
  \text{LEM} = 0.17423.
\end{equation*}
Thus, five stages give a smaller LEM than four stages.
\item Relaxing the nonnegativity constraint yields the optimized coefficients in Table~\ref{coeffs:optim4neg} with four stages.
\begin{table}[ht!]
\begin{center}
\begin{tabular}{|c||c|c|}
\hline
 & $a_{1,*}$ & $a_{2,*}$\\
\hline
1 & 0.39439914 & $-$0.092758759\\
\hline
2 & $-$0.12415477 & 0.60150021\\
\hline
3 & $-$0.10830436 & 0.13987854\\
\hline
4 & 0.83806 & 0.35138\\
\hline
\hline
 & $a_{3,*}$ & $a_{4,*}$\\
\hline
1 & 0.33190506 & 0.19579292\\
\hline
2 & 0.064464935 & 0.68067707\\
\hline
3 & 0.30064 & $-$0.13758\\
\hline
4 & 0.30299 & 0.26111\\
\hline
\end{tabular}
\caption{A four-operator splitting with some negative coefficients.\label{coeffs:optim4neg}}
\end{center}
\end{table}

The analysis shows
\begin{equation*}
  \text{LEM} = 0.80685.
\end{equation*}
The LEM is significantly smaller. If it should turn out to be uncritical to allow
slightly negative coefficients in the propagation of one of the sub-operators in (\ref{mhdeqns1})
without inducing instability, this method might significantly improve the accuracy and efficiency.
\end{itemize}

\section{Splitting into three operators}\label{subsec:3ops}

As explained earlier, a splitting into three operators can also be
taken into consideration, so in this section we give a method with
optimized LEM and a Milne pair for error estimation.

\begin{itemize}
\item An optimized three-operator splitting with positive coefficients is given in Table~\ref{coeffs:optim3pos}.
\begin{table}[ht!]
\begin{center}
\begin{tabular}{|c||c|c|c|}
\hline
 & $a_{1,*}$ & $a_{2,*}$ & $a_{3,*}$\\
\hline
1 & 0.31162504 & 0.27879542 & 0.67306805\\
2 & 2.4409272E-8 & 0.44755292 & 0.053280272\\
3 & 0.68837493 & 0.27365165 & 0.27365167\\
\hline
\end{tabular}
\caption{A three-operator splitting with positive coefficients.\label{coeffs:optim3pos}}
\end{center}
\end{table}
The local error measure is\nobreak
\begin{equation*}
\text{LEM} = 0.29596
\end{equation*}

\item Finally, here we present a Milne pair for asymptotically correct error estimation.
This construction requires eight compositions of the flow. The basic integrator is the
method \texttt{AK 3-2 (i)} from \cite{splithp} whose LEM equals 1.1. The adjoined method
is given in Table~\ref{tab:milne3ops}.
The numerically determined value for $\gamma$ in (\ref{eq:le-milne}) is $$\gamma=\frac{1}{4.1092266}.$$

\begin{table}[ht!]
\begin{center}
\begin{tabular}{|c||c|c|c|}
\hline
 & $a_{1,*}$ & $a_{2,*}$ & $a_{3,*}$\\
\hline
1 & 0.31133359 & 0.18034427 & 0.42236491\\
\hline
2 & 0 & 0.30701733 & 0.064464935\\
\hline
3 & 0.064046873 & 0.16776137 & 0.072032368\\
\hline
4 & 0 & 0.0046346054 & 0.10036361\\
\hline
5 & 0.23378298 & 0.043613842 & 0.043613842\\
\hline
6 & 0.18917015 & 0.00134284 & 0.14798253\\
\hline
7 & 0 & 0.19523257 & 0.048592873\\
\hline
8 & 0.20166638 & 0.10005315 & 0.10005315\\
\hline
\end{tabular}
\caption{A method with positive coefficients for three operators
adjoined to \texttt{AK 3-2 (i)} to form a Milne pair.\label{tab:milne3ops}}
\end{center}
\end{table}

The local error measure for the adjoined method is
\begin{equation*}
\text{LEM} = 0.12167.
\end{equation*}

\end{itemize}

\section{Numerical example -- the viscous Burgers equation}\label{sec:burgers}

The simulation models arising in astrophysical applications are extremely complex
and equally challenging in a theoretical error analysis. For numerical tests
which serve as a proof of concept for the adaptive splitting approach,
we therefore resort to a simpler parabolic--hyperbolic problem which can
readily be solved by adaptive two-operator splitting.

A representative test example which is a simplified
version of the Navier--Stokes equation and thus catches essential features
of incompressible flow is given by the viscous Burgers' equation \cite{evans98}
\begin{equation}\label{burgers1}
\partial_t u(x,t) = \kappa\partial_{xx}u(x,t) - \cdot u(x,t)\partial_x u(x,t), \quad u(x,0)=u_0(x),
\end{equation}
where we chose $\kappa=0.01.$
This has the known time-independent exact solution $u(x,t)=-\tanh(x)$ in the
\emph{shock case} and a known time-dependent solution in the \emph{rarefaction case},
see \cite{holdenetal10}. The question of the convergence of time semi-discretization
in appropriate Banach spaces has been addressed in \cite{holdenetal13}.

In the simplest setting, the problem is solved on a bounded interval with periodic boundary conditions.
We have implemented second-order splittings for this case to demonstrate that adaptive
time-stepping correctly reflects the solution behavior. The parabolic part defined by
$A=\partial_{xx}$ is discretized in space by a Fourier spectral method with 4096 modes,
and the nonlinear operator $B$ by the Lax--Wendroff finite difference scheme.
The splittings used constitute a Milne pair of order two specified at
\cite{splithp} with denotation \texttt{Symm-Milne-32}.
Between substeps, the fast Fourier transform serves to switch between real space and frequency space.

Table~\ref{table:order} shows the empirical convergence order of the global error at $t=0.28174$, computed with respect to
a precise numerical reference solution, for both the splitting methods from the Milne pair,
where we have chosen $u_0(x)=\tfrac12 \ee^{1/(x^2-1)},\ x\in (-1,1)$.
This serves to verify that the methods indeed show their theoretical orders in our setting
and can be used for reliable local error estimation. The error has been computed
as compared to a reference solution computed with stepsize $1/100$ times the smallest step, i.e. $9.78\cdot 10^{-6}$.

\begin{table}
\begin{center}
{\tiny
\begin{tabular}{|l||r|r|r||r|r|r|}
\hline
$h$ & global error Milne22 & order & error constant & global error Milne32 & order & error constant\\ \hline
0.0625 & 2.14e-08 & 2.30 & 1.28e-05 & 8.10e-09 & 1.70 & 9.01e-07\\
0.03125 & 6.29e-09 & 2.20 & 1.27e-05 & 3.01e-09 & 1.87 & 1.97e-06\\
0.01563 & 1.61e-09 & 2.24 & 1.59e-05 & 7.71e-10 & 1.98 & 2.95e-06\\
0.007813 & 2.89e-10 & 2.06 & 1.49e-05 & 1.34e-10 & 1.96 & 1.83e-06\\
0.003906 & 4.95e-11 & 2.00 & 4.46e-06 & 2.76e-11 & 1.94 & 1.28e-06\\
0.001953 & 1.30e-11 & 1.98 & 3.46e-06 & 7.72e-12 & 1.97 & 1.72e-06\\
0.0009766 & 2.72e-12 & 1.93 & 2.54e-06 & 1.65e-12 & 1.97 & 1.36e-06 \\
\hline
\end{tabular}
}
\caption{Order of the global error of second order Milne splittings applied to (\ref{burgers1}).\label{table:order}}
\end{center}
\end{table}

\begin{figure}
\begin{center}
\includegraphics[width=1.0\textwidth]{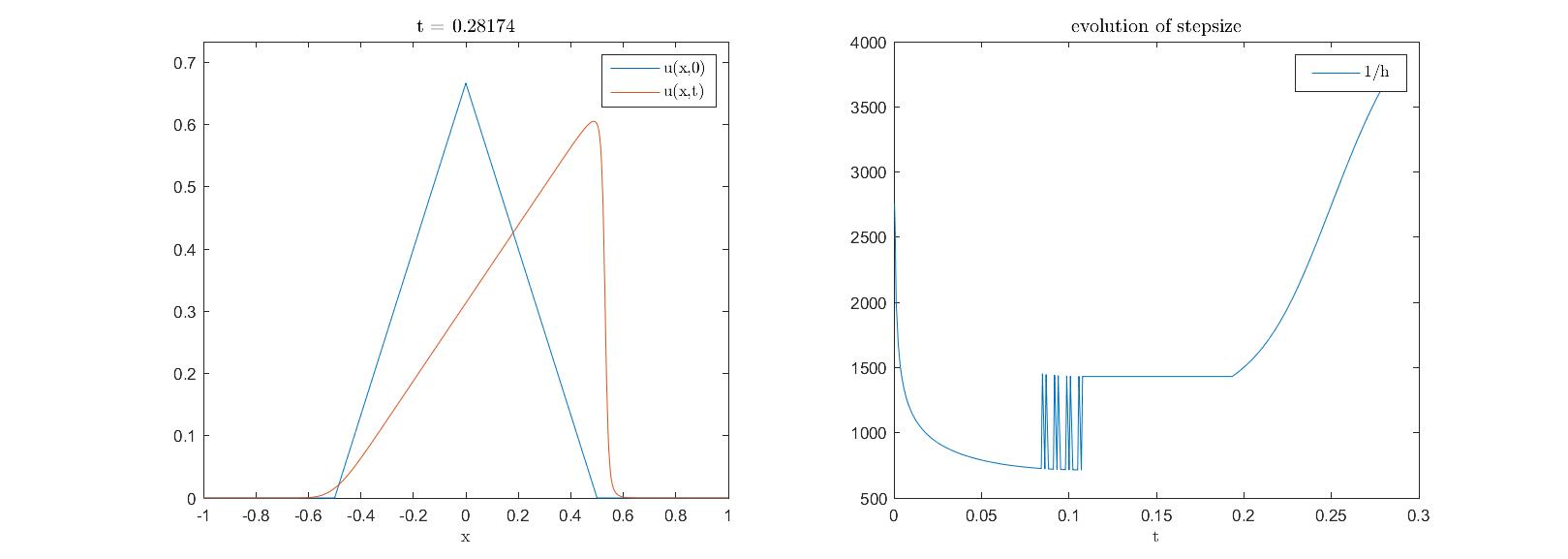}
\caption{Solutions $u(x,0)$ and $u(x,0.28174)$ for (\ref{burgers1}) (left)
and reciprocal of adaptively generated stepsizes (right).\label{figure:stepsizes}}
\end{center}
\end{figure}

On this basis, we solved (\ref{burgers1}) with adaptive step-size control with
a prescribed tolerance of $10^{-5}$. For this purpose, we used as initial condition
a symmetrical hat function of height 3/2, centered at 0. For this datum,
a steep gradient in the solution develops, whence the consistent behavior of the adaptive time-stepper
can be validated. In Figure~\ref{figure:stepsizes}, $u_0(0)$ and $u(x,t)$ at
$t=0.28174$ are shown in the left plot, and in the right plot
$1/h$ for the adaptively chosen stepsizes $h$. These grow from a cautious initial
guess while the solution is smooth, and after transient oscillations decrease to the prescribed lower limit
when the solution features a shock-like behavior.

\section{Conclusions and Outlook}

We have put forward adaptive multi-operator splitting methods based on
asymptotically correct estimation of the local time-stepping error. We have constructed optimized
methods and pairs of methods for error estimation satisfying a
physically mandatory positivity condition. To demonstrate that an adaptive
strategy for time-stepping correctly captures the possible loss
of smoothness for parabolic-hyperbolic flow problems by consistently
reducing step-sizes, we solved the viscous Burgers equation.
Future work will comprise the construction of Milne pairs also for
four operator splitting. The construction of methods for splitting
into more operators, which may provide advantages for MHD simulations,
is limited by computational resources due to the exponential growth
of the number of order conditions.

\section*{Acknowledgements}

This work was supported by the Austrian Science Fund (FWF) under the grants P 35485-N
(O. Koch, F. Kupka) and P 33140-N (F. Kupka).

We would like to thank J. Deimel and L. H{\"o}ltinger,
students at TU Wien, for contributions to the construction of methods.

O. Koch and F. Kupka gratefully acknowledge being granted the status of
Research Fellow at the Faculty of Mathematics, University of Vienna.


\end{document}